\documentclass[conference,a4paper]{IEEEtran}
\IEEEoverridecommandlockouts

\usepackage{amsmath,amssymb,amscd}

\def\l{\lambda}

\def\e{\epsilon}
\def\o{\omega}

\def\r{\rho}

\def\bp{\begin{proposition}}
\def\ep{\end{proposition}}
\def\bt{\begin{theo}}
\def\et{\end{theo}}
\def\be{\begin{equation}}
\def\ee{\end{equation}}
\def\bl{\begin{lemma}}
\def\el{\end{lemma}}
\def\bc{\begin{corollary}}
\def\ec{\end{corollary}}
\def\pr{\noindent{\bf Proof: }}

\def\bd{\begin{definition}}
\def\ed{\end{definition}}

\def\l{\ell}

\newtheorem{theo}{Theorem}[section]
\newtheorem{lemma}{Lemma}[section]
\newtheorem{definition}{Definition}[section]
\newtheorem{corollary}{Corollary}[section]
\newtheorem{proposition}{Proposition}[section]

\begin{document}

\title{Decoupling of Fourier Reconstruction System for Shifts of Several Signals}

\author{\IEEEauthorblockN{Yosef Yomdin\IEEEauthorrefmark{1}\IEEEauthorrefmark{2},
Niv Sarig\IEEEauthorrefmark{1}\IEEEauthorrefmark{3} and
Dmitry Batenkov\IEEEauthorrefmark{1}\IEEEauthorrefmark{4}}
\IEEEauthorblockA{\IEEEauthorrefmark{1}Department of Mathematics,Weizmann Institute of Science, Rehovot 76100, Israel}
\IEEEauthorblockA{\IEEEauthorrefmark{2}Email: yosef.yomdin@weizmann.ac.il}
\IEEEauthorblockA{\IEEEauthorrefmark{3}Email: niv.sarig@weizmann.ac.il}
\IEEEauthorblockA{\IEEEauthorrefmark{4}Email: dima.batenkov@weizmann.ac.il}
\thanks{This research was supported by the Adams Fellowship Program of the Israeli Academy of Sciences and Humanities, 
ISF Grant No. 639/09, and by the Minerva foundation. We would like to thank the referees for useful corrections and remarks.}}

\maketitle

\begin{abstract}
We consider the problem of ``algebraic reconstruction'' of linear combinations of shifts of several signals
$f_1,\ldots,f_k$ from the Fourier samples. For each $r=1,\ldots,k$ we choose sampling set $S_r$ to be a subset of the
common set of zeroes of the Fourier transforms ${\cal F}(f_\l), \ \l \ne r$, on which ${\cal F}(f_r)\ne 0$. We show that
in this way the reconstruction system is reduced to $k$ separate systems, each including only one of the signals $f_r$.
Each of the resulting systems is of a ``generalized Prony'' form. We discuss the problem of unique solvability of such
systems, and provide some examples.

\end{abstract}


\section{Introduction}
\setcounter{equation}{0}

In this paper we consider reconstruction of signals of the following a priori known form:

\be\label{equation_decoupling_model}
F(x)=\sum_{j=1}^k \sum_{q=1}^{q_j} a_{jq} f_j(x-x_{jq}),
\ee
with $a_{jq}\in \mathbb{R}, \ x_{jq}=(x^1_{jq},\ldots,x^n_{jq}) \in {\mathbb R}^n.$
We assume that the signals $f_1,\dots, f_k: \mathbb{R}^n\to\mathbb{R}$ are known (in particular, their Fourier transforms ${\cal F}(f_j)$ are known),
while $a_{jq}, \ x_{jq}$ are the unknown signal parameters, which we want to find from Fourier samples of $F$. We explicitly
assume here that $k\geq 2$. So the usual methods which allow one to solve this problem ``in closed form'' in the case of shifts
of a single function (see \cite{Vet,Bat.Sar.Yom,Sar}) are not directly applicable. Still, we shall show that in many cases an
explicit reconstruction from a relatively small collection of Fourier samples of $F$ is possible. Practical importance of signals
as above is well recognized in the literature: for some discussions and similar settings see, e.g. 
\cite{Vet,gedalyahu2011multichannel, peter2011nonlinear}.

We follow a general line of the ``Algebraic Sampling'' approach (see \cite{Vet,sig_ack,Bat.Yom1} and references therein), i.e. we
reconstruct the values of the unknown parameters, solving a system of non-linear equations, imposed by the measurements (system
(\ref{equation_uncoupled_system}) below). The equations in this system appear as we equate the ``symbolic'' expressions of the
Fourier samples, obtained from (\ref{equation_decoupling_model}), to their actual measured values.

Our specific strategy is as follows: we choose a sampling set $S_r \subset {\mathbb R}^n, \  r=1,\ldots,k,$ in a special way,
in order to ``decouple'' (\ref{equation_uncoupled_system}), and to reduce it to $k$ separate systems, each including only one
of the signals $f_r$. To achieve this goal we take $S_r$ to be a subset of the common set of zeroes of the Fourier transforms
${\cal F}(f_\l), \ \l\ne r$.

The decoupled systems turn out to be of a ``generalized Prony'' type:
\be\label{equation_decoupled1}
\sum_{j=1}^{N} a_j y_j^{s_\l} =  m_\l, \quad \l=1,2,\dots, \ s_\l \in S\subset {\mathbb R}^n.
\ee
The standard Prony system, where the sample set $S$ is the set of integer points in a cube of a prescribed size, allows for a 
solution ``in closed form'' (see, for example, \cite{Bat.Sar.Yom,rao1992mbp,Sar,stoica2005spectral} and references therein). We 
are not aware of any method for an explicit solution of generalized Prony systems. However, ``generic'' solution methods can be
applied. Their robustness can be estimated via Tur\'an-Nazarov inequality for exponential polynomials and its discrete version 
(\cite{Fri.Yom,Naz}). Some initial results in this direction have been presented in \cite{Sar,Bat.Sar.Yom}. Below we further 
extend these results, restricting ourselves to the uniqueness problem only.

\section{Reconstruction System and its Decoupling}\label{subsection_decoupling_system}
\setcounter{equation}{0}

For $F$ of the form (\ref{equation_decoupling_model}) and for any $s=(s^1,\ldots,s^n)\in {\mathbb R}^n$ we have for the
sample of the Fourier transform ${\cal F}(F)$ at $s$

\begin{eqnarray*}
&{\cal F}(F)(s)&= \int_{{\mathbb R}^n} e^{-2\pi isx} F(x)dx\\
&&=\sum_{j=1}^k \sum_{q=1}^{q_j} a_{jq} e^{-2\pi isx_{jq}}{\cal F}(f_j)(s).
\end{eqnarray*}
So taking samples at the points $s_\l=(s^1_\l,\ldots,s^n_\l)$ of the sample set $S=\{s_1,\dots,s_m\}$, and denoting the vector
$e^{-2\pi ix_{jq}}=(e^{-2\pi ix^1_{jq}},\dots,e^{-2\pi ix^n_{jq}})$ by $y_{jq}=(y^1_{jq},\dots, y^n_{jq})$ we get our
reconstruction system in the form

\be\label{equation_uncoupled_system}
\sum_{j=1}^k \sum_{q=1}^{q_j} a_{jq} {\cal F}(f_j)(s_\l)y_{jq}^{s_\l} = {\cal F}(F)(s_\l), \ \l=1,\dots,m,
\ee
in the standard multi-index notations. In system (\ref{equation_uncoupled_system}) the right hand sides ${\cal F}(F)(s_\l)$
are the known measurements, while the Fourier samples ${\cal F}(f_j)(s_\l)$ are known by our assumptions. The unknowns in
(\ref{equation_uncoupled_system}) are the amplitudes $a_{jq}$ and the shifts $x_{jq}$, encoded in the vectors $y_{jq}$.

In the case $k=1$ we could divide the equations in (\ref{equation_uncoupled_system}) by ${\cal F}(f_1)(s_\l)$ and obtain
directly a Prony-like system. However, for $k\geq 2$ this transformation usually is not applicable. Instead we ``decouple" system (\ref{equation_uncoupled_system})
with respect to the signals $f_1,\ldots,f_k$ using the freedom in the choice of the sample set $S$. Let
\[Z_\l=\bigl\{x\in {\mathbb R}^n, \ {\cal F}(f_\l)(x)=0\bigr\}\] denote the set of zeroes of the Fourier transform ${\cal F}(f_\l)$.
For each $r=1,\dots,k$ we take the sampling set $S_r$ to be a subset of the set \[W_r=(\bigcap_{\l\ne r} Z_\l)\setminus Z_r\]
of common zeroes of the Fourier transforms ${\cal F}(f_\l), \ \l\ne r$, but not of ${\cal F}(f_r)$. For such $S_r$ all the
equations in (\ref{equation_uncoupled_system}) vanish, besides those with $j=r$. Hence we obtain:

\bp\label{proposition_coupling}
Let for each $r=1,\dots,k$ the sampling set $S_r$ satisfy \[S_r=\{s_{r1},\dots,s_{rm_r}\}\subset W_r.\] Then for each
$r$ the corresponding system (\ref{equation_uncoupled_system}) on the sample set $S_r$ takes the form

\be\label{equation_decoupled}
\sum_{q=1}^{q_r} a_{rq} y_{rq}^{s_{r\l}} =  c_{r\l}(F), \ \l=1,\dots,m_r,
\ee
where $c_{r\l}(F)= {{{\cal F}(F)(s_{r\l})}/{{\cal F}(f_r)(s_{r\l})}}$. $\square$
\ep

\smallskip

So (\ref{equation_uncoupled_system}) is decoupled into $k$ generalized Prony systems (\ref{equation_decoupled}), each
relating to the shifts of the only signal $f_r$. The problem is that some (or all) of the sets $W_r$ may be too small,
and the resulting systems (\ref{equation_decoupled}) will not allow us to reconstruct the unknowns $a_{rq}$ and $y_{rq}$.
Another problem is instability of zero finding, which may lead to only approximate zeroes of Fourier transforms. We have
at present only initial results outlying applicability of the Fourier decoupling method (\cite{Sar}). In a ``good'' case 
where the zero sets $Z_\l$ of the Fourier transforms ${\cal F}(f_\l), \ \l=1,\ldots,k,$ are nonempty $n-1$-dimensional 
hypersurfaces meeting one another transversally, still for $k>n+1$ the intersection of $Z_\l, \ \l\ne r,$ is empty. So the 
resulting systems (\ref{equation_decoupled}) contain no equations. Hence we can apply the above decoupling only for $k \leq n+1$. 

\medskip

Some specific examples, as well as investigation of the conditions on $f_1,\ldots,f_k$ which provide solvability of systems
(\ref{equation_decoupled}) were presented in \cite{Sar}. In one-dimensional case $(n=1, k=2)$ these conditions can be given
explicitly. In this case $W_1=W_1(f_1,f_2)$ consists of zeroes of ${\cal F}(f_2)$ which are not zeroes of ${\cal F}(f_1)$,
and $W_2=W_2(f_1,f_2)$ consists of zeroes of ${\cal F}(f_1)$ which are not zeroes of ${\cal F}(f_2)$. The following result
has been proved (for real Prony systems) in \cite{Sar}. Here we extend it to the case of system (\ref{equation_decoupled})
which has purely imaginary exponents. The constant $2N$ below is sharp, in contrast with the constant $C(n,d)$ in 
(multidimensional) Theorem \ref{span} below.

\smallskip

Let in (\ref{equation_decoupling_model}) \ $n=1,k=2,$ and let $q_1=q_2=N$. Assume that for the signals $f_1,f_2$ in
(\ref{equation_decoupling_model}) each of the sets $W_1$ and $W_2$ contains at least $2N$ elements. Let $D_j, \ j=1,2,$ be
the length of the shortest interval $\Delta_j$ such that $S_j=\Delta_j\cap W_j$ contains exactly $2N$ elements, and let
$\rho_j={1\over {D_j}}$.

\bt \label{solv.cond}
For shifts $x_{jq}$ in the interval $[0,\rho_j), \ j=1,2,$ systems (\ref{equation_decoupled}) with the sampling sets $S_1,S_2$
are uniquely solvable.
\et
\pr
Let us fix $j=1$. The proof for $j=2$ is the same. Substituting $y_{1q}=e^{-2\pi i x_{1q}}$ associates to a solution
$(a_{1q},y_{1q}), \ q=1,\ldots,N,$ of (\ref{equation_decoupled}) an exponential polynomial
$H(s)=\sum_{q=1}^N a_{1q}e^{-2\pi i x_{1q}s}$ with purely imaginary exponents. If (\ref{equation_decoupled}) has two
different solutions, the corresponding exponential polynomials $H_1(s)$ and $H_2(s)$ are equal for each $s \in S_1.$ Hence
$S_1$ is a set of zeroes of $H_2(s)-H_1(s)$, which is an exponential polynomial of the order at most $2N$. On the other hand, by
Langer's lemma (Lemma 1.3 in \cite{Naz}) such polynomial can have in each interval of length $D$ at most $2N-1+{{\rho D}\over {2\pi}}$
zeroes, where $\rho$ is the maximum of the absolute values of the exponents. In our case $D=D_1$ and
$\rho < 2\pi\rho_1={{2\pi}\over {D_j}}$. Hence ${{\rho D}\over {2\pi}}$ is strictly less than $1$, and so
the number of zeroes of $H_2-H_1$ is at most $2N-1$, in contradiction with the assumptions. $\square$

\section{Examples}\label{examples}
\setcounter{equation}{0}

Some examples of Fourier decoupling have been presented in \cite{Sar}. In these examples the sets $W_r$ are ``large enough''
to reduce the problem (with the number of allowed shifts fixed but arbitrarily large) to a set of decoupled standard Prony
systems.

\smallskip

In dimension one we can take, for example, $f_1$ to be the characteristic function of the interval $[-1,1],$ while
$f_2(x)=\delta(x-1)+\delta(x+1).$ So we consider signals of the form

\be\label{equation_decoupling_model1}
F(x)= \sum_{q=1}^{N} [a_{1q} f_1(x-x_{1q})+a_{2q} f_2(x-x_{2q})].
\ee
Easy computations show that

\[
{\cal F}(f_1)(s)=\sqrt{\frac{2}{\pi}}\frac{\sin s}{s}
\]
and
\[
{\cal F}(f_2)(s)=\sqrt\frac{2}{\pi}\cos s.
\]
So the zeros of the Fourier transform of $f_1$ are the points $\pi n,\ n\in {\mathbb Z}\setminus \{0\}$ and those of $f_2$
are the points $({1\over 2}+n)\pi,\ n\in {\mathbb Z}$. These sets do not intersect, so $W_1=\{\pi n\}$, and
$W_2=\{({1\over 2}+n)\pi\}$. Since $W_1$ and $W_2$ are just shifted integers ${\mathbb Z}$, the generalized Prony systems
in (\ref{equation_decoupled}) are actually the standard ones. For $f_2$ the system (\ref{equation_decoupled}) takes the form
\[
\frac{{\cal F}(F)(\pi n)}{\sqrt\frac{2}{\pi}(-1)^n}=\sum_{q=1}^Na_{2q}({y_{2q}})^{ \pi n}, \ n\in {\mathbb Z}.
\]
If we denote $M_n=\frac{{\cal F}(F)(\pi n)}{\sqrt\frac{2}{\pi}(-1)^n}$ , $A_q=a_{2q}(y_{2q})^\pi$
and $\eta_q=(y_{2q})^\pi$ we get the usual Prony system
\[
M_n=\sum_{q=0}^NA_q\eta_q^n\ , \ n\in {\mathbb Z}.
\]
For $f_1$ we get
\[
\frac{{\cal F}(F)(({1\over 2} +n)\pi)}{\sqrt\frac{2}{\pi}\frac{(-1)^{n+1}}
{({1\over 2} +n)\pi}}=\sum_{q=1}^Na_{1q}({y_{1q}})^{({1\over 2} +n)\pi}  \ , \ n\in {\mathbb Z}.
\]
In this case we denote $\mu_n=\frac{{\cal F}(F)(({1\over 2} +n)\pi)}{\sqrt\frac{2}{\pi}
\frac{(-1)^{n+1}}{({1\over 2} +n)\pi}},\ \alpha_q=a_{1q}(y_{1q})^{\frac\pi2}$ and
$\xi_q=(y_{1q})^\pi$ and we get again the usual Prony system
\[
\mu_n=\sum_{q=1}^N\alpha_q\xi_q^n,\ n\in {\mathbb Z}.
\]
Solving these two systems by any standard method will give us the translations and amplitudes of the functions $f_1, f_2$.
Notice that a possible non-uniqueness of the solutions is imposed here by the substitutions $\eta_q=(y_{2q})^\pi$ and
$\xi_q=(y_{1q})^\pi$.

\smallskip

In dimension two we can take, in particular, $f_1,f_2,f_3$ to be the characteristic functions of the three
squares: $Q_1=[-3,3]^2, \ Q_2=[-5,5]^2,$ and $Q_3$ which is the rotation of the square $[-\sqrt 2,\sqrt 2]^2$ by
${\pi \over 4}$. So we put

\be
\chi_j(x) =\left\{\begin{array}{lr}1&x\in Q_j
\\0&x\not\in Q_j
\end{array}\right.
\ee
and consider signals of the form

\be\label{equation_decoupling_model2}
F(x)=\sum_{j=1}^3 \sum_{q=1}^{q_j} a_{jq} \chi_j(x-x_{jq}), \quad  \text{ with } a_{jq}\in \mathbb{R},\; x_{jq} \in {\mathbb R}^3.
\ee
The following result is proved in \cite{Sar}:

\smallskip

\bp
The zero sets  $Z_1,Z_2$ and $Z_3$ of the Fourier transforms of the three functions $\chi_1,\chi_2$ and $\chi_3$ intersect each
other in such a way that the decoupling procedure based on the sets $W_1=(Z_2\cap Z_3)\setminus Z_1, W_2=(Z_3\cap Z_1)\setminus Z_2$
and $W_3=(Z_1\cap Z_2)\setminus Z_3$ provides three standard Prony systems for the
shifts of each of the functions.
\ep

\smallskip

\noindent{\bf Sketch of the proof:} Simple calculation gives
\be\begin{array}c
{\cal F}(\chi_1)(\o,\r)=4\frac{\sin 3\o}{\o}\cdot\frac{\sin 3\r}{\r}\\
{\cal F}(\chi_2)(\o,\r)=4\frac{\sin 5\o}{\o}\cdot\frac{\sin 5\r}{\r}\\
{\cal F}(\chi_3)(\o,\r)=8\frac{\sin \frac{\o+\r}2}{\frac{\o+\r}2}\cdot\frac{\sin \frac{\o-\r}2}{\frac{\o-\r}2}.
\end{array}
\ee
So $Z_1$ is the union of horizontal or vertical lines crossing the Fourier plane's axes at $(0,\frac{n\pi}{3})$ or
$(\frac{n\pi}{3},0)$ respectively, for all non zero integer $n$. Similarly for $Z_2$, with the only difference that the lines cross
the axes at $(0,\frac{n\pi}5)$ or $(\frac{n\pi}5,0)$. \\$Z_3$ is the union of lines with slopes $1$ or $-1$ crossing the $\o$ axis at
$2\pi n$ for some non zero integer $n$. Hence for any two integers $n$ and $m$ we have
$(\frac{1+5n}5,\frac{1+5n}5)\in S_1, (\frac{1+3m}3,\frac{1+3m}3)\in S_2$ and since $\frac{1+3m}3\pm\frac{1+5n}5$ is not an integer,
$(\frac{1+3m}3,\frac{1+5n}5)\in S_3$. These three points form a triangle which repeats itself as a periodic pattern. Appropriate
transformations now bring the decoupled systems (\ref{equation_decoupled}) to the form of the standard two-dimensional Prony system.
See \cite{Sar,Bat.Sar.Yom} for a new approach to solving such systems and for the results of numerical simulations. $\square$

\section{Uniqueness of Reconstruction}\label{Non-uniform}
\setcounter{equation}{0}

Application of Proposition \ref{proposition_coupling} prescribes the choice of sample points from the common zeroes of
the Fourier transforms ${\cal F}(f_j)$. So the geometry of the sample sets $S_r$ may be complicated, and the known results
on unique solvability of the standard Prony system (\cite{Bat.Sar.Yom,Bat.Yom,rao1992mbp,stoica2005spectral}) are not directly
applicable. Non-Uniform Sampling in Prony-type systems is also essential in other problems of algebraic signal reconstruction.
In particular, recently it appeared as a key point in a proof of the Eckhoff conjecture, related to the accuracy of
reconstruction of piecewise-smooth functions from their Fourier samples (\cite{Bat}).

There are results on a behavior of exponential polynomials on arbitrary sets, which can provide important information on unique
solvability and robustness of the generalized Prony system. In particular, this concerns the Turan-Nazarov inequality
(\cite{Naz}), and its extension to discrete sets obtained in \cite{Fri.Yom}. In this last paper for each set $S$ a quantity
$\o_D(S)$ has been introduced, measuring, essentially, the robustness of solvability of a generalized Prony system with the
sample points $s_\l\in S$. Here $D$ comprises the ``discrete'' parameters of the Prony system to be solved. $\o_D(S)$ can be
explicitly estimated in terms of the metric entropy of $S$ (see below), and we expect that in many important cases the quantity
$\o_D(W_r)$ for the zeroes sets $W_r$ of the Fourier transforms ${\cal F}(f_j)$ can be effectively bounded from below. Some initial
results and discussions in this direction, mainly in dimension one, are presented in \cite{Sar,Bat.Yom1}. In the present paper we
do not consider robustness of the Prony system, but provide a new multi-dimensional result on the uniqueness of solutions, in the
lines of \cite{Sar,Fri.Yom} and Theorem \ref{solv.cond} above.

\smallskip

Let us recall that for $Z$ a bounded subset of ${\mathbb R}^n,$ and for $\e>0$ the covering number $M(\e,Z)$ is the minimal number
of $\e$-balls in ${\mathbb R}^n,$ covering $Z$. The $\e$-entropy $H(\e,Z)$ is the binary logarithm of $M(\e,Z)$.

\smallskip

Let $H(s)=\sum_{j=1}^d a_je^{\lambda_j\cdot s},$ with
$a_j\in {\mathbb R}, \ \lambda_j=(\lambda_{j1},\ldots,\lambda_{jn})\in {\mathbb R}^n,$ be a real exponential polynomial in
$s\in {\mathbb R}^n.$ Denote $Z(H)$ the set of zeroes of $H$ in ${\mathbb R}^n$, and let $Q^n_R$ be the cube in ${\mathbb R}^n$ with
the edge $R$. The following result is a special case of Lemma 3.3 proved in \cite{Fri.Yom}:

\bp \label{entropy}
For each $R>0,$ and $\e$ with $R>\e>0$ we have $M(\e,Z(H)\cap Q^n_R)\leq C(d,n)({R\over \e})^{n-1}$. $\square$
\ep

\smallskip

The explicit expression for $C(d,n)$ is given in \cite{Fri.Yom}, via Khovanski's bound (\cite{Kho}) for ``fewnomial'' systems.
Consider now a generalized Prony system (\ref{equation_decoupled1}) with a finite set $S$ of samples allowed:

\be \label{Pron}
\sum_{j=1}^{N} a_j y_j^{s_\l} =  m_\l, \ s_\l \in S=\{s_1,\ldots,s_m \}\subset {\mathbb R}^n.
\ee
We shall consider only real solutions of (\ref{Pron}) with $y_j$ having all its coordinates positive.

\bt \label{span}
Let $S=\{s_1,\ldots,s_m \}\subset Q_R^n$ be given, such that for a certain $\e>0$ we have $M(\e,S)> C(2N,n)({R\over \e})^{n-1}.$
Then system (\ref{Pron}) has at most one solution.
\et
\pr
Associate to a solution $(a_j,y_j), \ j=1,\ldots,N,$ of (\ref{Pron}) an exponential polynomial
$H(s)=\sum_{j=1}^N a_je^{\lambda_j\cdot s},$ where $y_j=e^{\lambda_j}, \ \lambda_j\in {\mathbb R}^n.$ If (\ref{Pron}) has two
different solutions, the corresponding exponential polynomials $H_1(s)$ and $H_2(s)$ are equal for each $s=s_\l\in S.$ Hence
$S$ is a set of zeroes of $H_2(s)-H_1(s)$, which is an exponential polynomial of order at most $2N$. By Proposition \ref{entropy}
we have $M(\e,S) \leq C(2N,n)({R\over \e})^{n-1}$ for each $\e>0$, in contradiction with the assumptions of the theorem. $\square$

\smallskip

Informally, Theorem \ref{span} claims that finite sets $S$ which cover (in a ``resolution $\e$'', for some $\e>0$), a significant
part of the cube $Q_R^n$, are uniqueness sets of the Prony system. The condition of Theorem \ref{span} on the sampling set $S$ is
quite robust with respect to the geometry of $S$, so we can explicitly verify it in many cases. In particular, for non-regular
lattices we get the following result:

\bd
For fixed positive $\alpha < {1\over 2}$ and $h>0,$ a set $Z' \subset {\mathbb R}^n$ is called an $(\alpha,h)$-net if it possesses
the following property: there exists a regular grid $Z$ with the step $h$ in ${\mathbb R}^n$ such that for each $z'\in Z'$ there
is $z\in Z$ with $||z'-z||\leq \alpha h,$ and for each $z\in Z$ there is $z'\in Z'$ with $||z'-z||\leq \alpha h.$
\ed
\bc \label{Nonreg}
Let $Z' \subset {\mathbb R}^n$ be an $(\alpha,h)$-net. Then for $R> C(2N)h(1-2\alpha)^{1-n}$ the set $S=Z\cap Q_R^n$ is a uniqueness
set of the Prony system \eqref{Pron}.
\ec
\pr
By definition, for each $z\in Z$ we can find $z'\in Z'$ inside the $\alpha h$-ball around $z$. Clearly, any two such points are
$h'=(1-2\alpha)h$-separated. So for each $\e < h'$ we have $M(\e,S)\geq |Z\cap Q_R^n|=({R\over h})^n.$ We conclude that
the inequality $({R\over h})^n > C(2N)({R\over h'})^{n-1},$ or $R > C(2N)h(1-2\alpha)^{1-n}$ implies the condition of Theorem
\ref{span}. $\square$

\smallskip

The condition of Theorem \ref{span} can be verified in many other situations, under natural assumptions on the sample set $S$.
In particular, using integral-geometric methods developed in \cite{Com.Yom}, it can be checked for the zero sets of Fourier
transforms of various types of signals. We plan to present these results separately.

\smallskip

\noindent{\bf Remark} The restriction to only positive solutions of Prony system is very essential for the result of Theorem
\ref{span}. Indeed, consider the Prony system

\be \label {Skolem}
a_1x_1^k+a_2x_2^k=m_k, \ k=0,1,\ldots.
\ee
If we put $a_1=1, \ x_1=1, \ a_2=-1, \ x_2=-1$, then $m_k=1^k-(-1)^k=0$ for each even $k$. So the regular grid of even integers
is not a uniqueness set for system (\ref{Skolem}). This fact is closely related to the classical Skolem-Mahler-Lech Theorem
(see \cite{Lec,Mey.vdP,Tao} and references therein) which says that the integer zeros of an exponential polynomial are the union of
complete arithmetic progressions and a finite number of exceptional zeros.  So such sets may be non-uniqueness sample sets for
complex Prony systems.

\smallskip

The proof of the Skolem-Mahler-Lech Theorem is relied on non-effective arithmetic considerations. Recently the problem of
obtaining effective such theorem was discussed in \cite{Tao}. This problem may turn to be important for understanding of
complex solutions of Prony systems. One can wonder whether the methods of Khovanskii (\cite{Kho}) and Nazarov (\cite{Naz}), as
well as their combination in \cite{Fri.Yom}, can be applied here.

\vskip1cm

\bibliographystyle{plain}

\end{document}